\documentclass{article}

\usepackage{amssymb,latexsym,amsmath}

\usepackage[pdftex]{graphicx}

\usepackage{hyperref}

\begin{document}

\pagestyle{plain}

\newtheorem{theorem}{Theorem}[section]

\newtheorem{proposition}[theorem]{Proposition}

\newtheorem{lema}[theorem]{Lemma}

\newtheorem{corollary}[theorem]{Corollary}

\newtheorem{definition}[theorem]{Definition}

\newtheorem{remark}[theorem]{Remark}

\newtheorem{exempl}{Example}[section]

\newenvironment{exemplu}{\begin{exempl}  \em}{\hfill $\square$

\end{exempl}}

\newcommand{\ea}{\mbox{{\bf a}}}

\newcommand{\eu}{\mbox{{\bf u}}}

\newcommand{\ueu}{\underline{\eu}}

\newcommand{\ueo}{\overline{u}}

\newcommand{\oeu}{\overline{\eu}}

\newcommand{\ew}{\mbox{{\bf w}}}

\newcommand{\ef}{\mbox{{\bf f}}}

\newcommand{\eF}{\mbox{{\bf F}}}

\newcommand{\eC}{\mbox{{\bf C}}}

\newcommand{\en}{\mbox{{\bf n}}}

\newcommand{\eT}{\mbox{{\bf T}}}

\newcommand{\eL}{\mbox{{\bf L}}}

\newcommand{\eR}{\mbox{{\bf R}}}

\newcommand{\eV}{\mbox{{\bf V}}}

\newcommand{\eU}{\mbox{{\bf U}}}

\newcommand{\ev}{\mbox{{\bf v}}}

\newcommand{\eve}{\mbox{{\bf e}}}

\newcommand{\uev}{\underline{\ev}}

\newcommand{\eY}{\mbox{{\bf Y}}}

\newcommand{\eK}{\mbox{{\bf K}}}

\newcommand{\eP}{\mbox{{\bf P}}}

\newcommand{\eS}{\mbox{{\bf S}}}

\newcommand{\eJ}{\mbox{{\bf J}}}

\newcommand{\eB}{\mbox{{\bf B}}}

\newcommand{\eH}{\mbox{{\bf H}}}

\newcommand{\leb}{\mathcal{ L}^{n}}

\newcommand{\eI}{\mathcal{ I}}

\newcommand{\eE}{\mathcal{ E}}

\newcommand{\hen}{\mathcal{H}^{n-1}}

\newcommand{\eBV}{\mbox{{\bf BV}}}

\newcommand{\eA}{\mbox{{\bf A}}}

\newcommand{\eSBV}{\mbox{{\bf SBV}}}

\newcommand{\eBD}{\mbox{{\bf BD}}}

\newcommand{\eSBD}{\mbox{{\bf SBD}}}

\newcommand{\ecs}{\mbox{{\bf X}}}

\newcommand{\eg}{\mbox{{\bf g}}}

\newcommand{\paromega}{\partial \Omega}

\newcommand{\gau}{\Gamma_{u}}

\newcommand{\gaf}{\Gamma_{f}}

\newcommand{\sig}{{\bf \sigma}}

\newcommand{\gac}{\Gamma_{\mbox{{\bf c}}}}

\newcommand{\deu}{\dot{\eu}}

\newcommand{\dueu}{\underline{\deu}}

\newcommand{\dev}{\dot{\ev}}

\newcommand{\duev}{\underline{\dev}}

\newcommand{\weak}{\stackrel{w}{\approx}}

\newcommand{\mild}{\stackrel{m}{\approx}}

\newcommand{\lrightarrow}{\stackrel{L}{\rightarrow}}

\newcommand{\rrightarrow}{\stackrel{R}{\rightarrow}}

\newcommand{\strong}{\stackrel{s}{\approx}}

\newcommand{\weakdown}{\rightharpoondown}

\newcommand{\opg}{\stackrel{\mathfrak{g}}{\cdot}}

\newcommand{\opunu}{\stackrel{1}{\cdot}}
\newcommand{\opdoi}{\stackrel{2}{\cdot}}

\newcommand{\opn}{\stackrel{\mathfrak{n}}{\cdot}}
\newcommand{\opx}{\stackrel{x}{\cdot}}

\newcommand{\tr}{\ \mbox{tr}}

\newcommand{\Ad}{\ \mbox{Ad}}

\newcommand{\ad}{\ \mbox{ad}}

\renewcommand{\contentsname}{ }

\title{Graphic lambda calculus and knot diagrams}

\author{Marius Buliga \\ 
\\
Institute of Mathematics, Romanian Academy \\
P.O. BOX 1-764, RO 014700\\
Bucure\c sti, Romania\\
{\footnotesize Marius.Buliga@imar.ro}}

\date{This version: 07.11.2012}

\maketitle

\begin{abstract}
In \cite{graphic} was proposed a graphic lambda calculus formalism, which has sectors corresponding to  untyped lambda calculus and emergent algebras. 

Here we explore the sector covering  knot diagrams, which are constructed as macros over the   graphic lambda calculus.
\end{abstract}

\section{Quick introduction and references}

The graphic lambda calculus \cite{graphic} is a formalism based on local or global  moves acting on locally planar trivalent graphs. In the mentioned paper we showed that "sectors" of this calculus are equivalent with untyped lambda calculus or with emergent algebras. (The formalism of emergent algebras \cite{buligairq} \cite{buligabraided} evolved from differential 
calculus on metric spaces with dilations \cite{buligadil1}.) 

For all the relevant notions and results consult \cite{graphic} for the graphic lambda calculus, \cite{lambdascale} for $\lambda$-Scale calculus (a first proposal of a calculus containing both untyped lambda calculus and emergent algebras). Those interested in metric spaces with dilations and their applications may consult the course notes \cite{cimpa} on sub-riemannian geometry from intrinsic viewpoint. For a larger view on the half-dreamed subject of "computing with space" see   \cite{buligachora}, where  a formalism  for emergent algebras based on decorated knot diagrams for emergent algebras was proposed. There we  called  "computing with space" the various manipulations of these diagrams with geometric content. (The interested reader may browse  various notes at \cite{post}, which are used as a repository and discussion place for the subject.) 

A significant part of the original contributions of this paper appeared as notes at  
\cite{graphicpost}.  

\paragraph{Acknowledgement.} This work was supported by a grant of the Romanian National Authority for Scientific Research, CNCS – UEFISCDI, project number 
PN-II-ID-PCE-2011-3-0383.

\section{Graphic lambda calculus} 

We introduce the primitives of graphic lambda calculus and then we start  from the last section of the paper \cite{graphic}. 

\paragraph{The set $GRAPH$.} It  consists of  oriented locally  planar graphs with decorated nodes, constructed  from a graphical alphabet of  elementary graphs, or "gates". 

\begin{definition}
The graphical alphabet contains the elementary graphs, or gates, denoted by $\lambda$, $\Upsilon$, $\curlywedge$, $\top$, and for any element $\varepsilon$ of a (given) commutative group $\Gamma$, a graph denoted by  $\bar{\varepsilon}$.  
\begin{enumerate}
\item[] $\lambda$ graph \hspace{1.cm} \includegraphics[width=20mm]{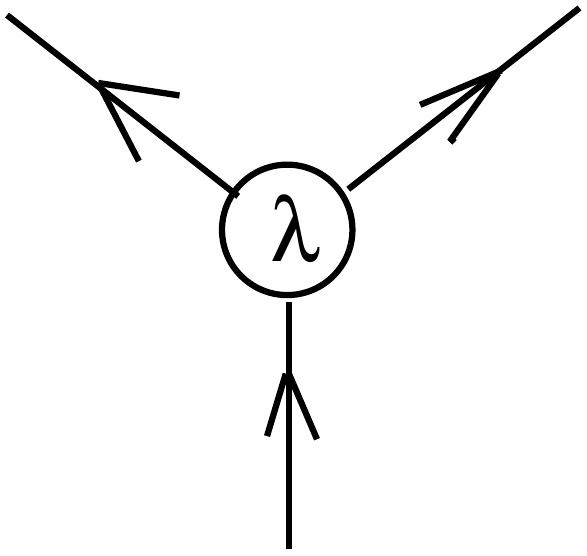}, \hspace{2.cm}
$\Upsilon$ graph \hspace{1.cm}\includegraphics[width=20mm]{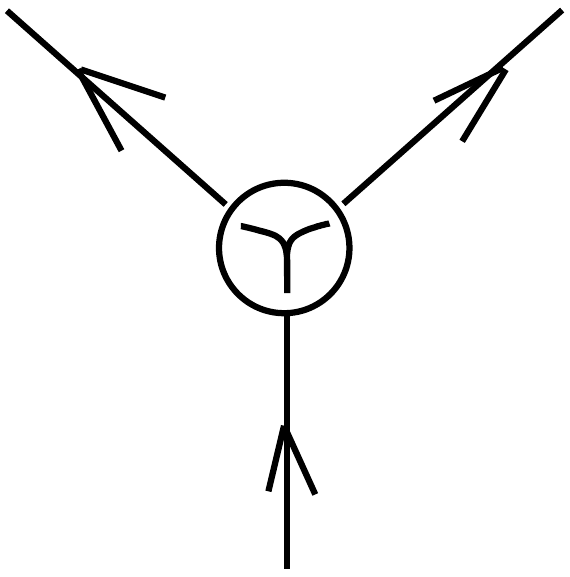}, 

\item[] $\curlywedge$ graph \hspace{1.cm}\includegraphics[width=15mm]{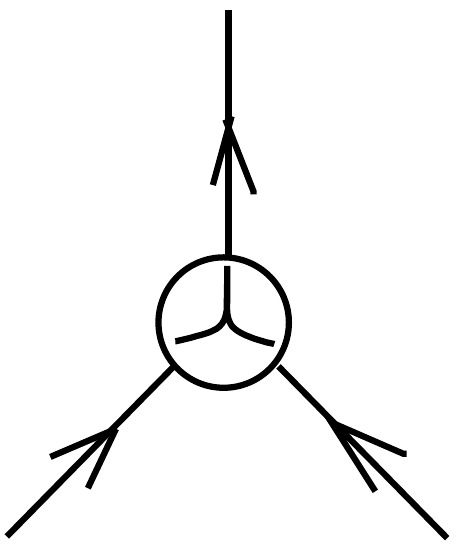}\hspace{1.cm}, \hspace{2.cm} $\bar{\varepsilon}$ graph \hspace{1.cm}\includegraphics[width=15mm]{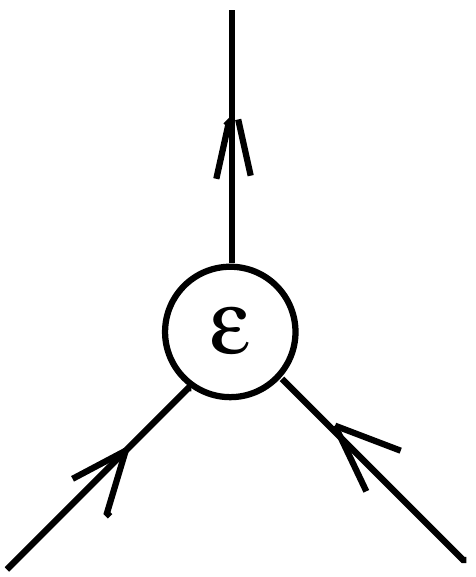}\hspace{.2cm}, 
\item[] $\top$ graph \hspace{1.cm}\includegraphics[width=8mm]{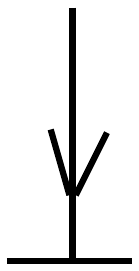}\hspace{1.cm}.
\end{enumerate}
With the exception of the $\top$, all other elementary graphs have three edges. The graph  $\top$ has only one edge. 
\label{defalp}
\end{definition}

$GRAPH$ is the set of oriented graphs obtained by grafting edges of a finite number of copies of the elements  of the graphical alphabet.For any node of the graph, the local embedding into the plane is given by the element of the graphical alphabet which decorates it. 

A graph $G \in GRAPH$ may have "free" edges, i.e. edges where no other elementary graph is grafted. We complete the graph $G \in GRAPH$ by adding to the free extremity of any free edge a decorated node, called "leaf",  with decoration "IN" or "OUT", depending on the orientation of the respective free edge. The set of leaves $L(G)$ decomposes into a disjoint union $L(G) = IN(G) \cup OUT(G)$ of in or out leaves.

By definition, a subgraph of a graph $G$ in $GRAPH$ is any subset of the reunion of nodes and edges of $G$. We don't  suppose that subgraphs are in  $GRAPH$, not even that a subgraph is a graph in the usual sense. For example,  a collection of edges of $G$, without any node, is a subgraph of $G$. 

We represent graphically a subgraph $P$ of a graph $G$ in $GRAPH$ in the following way. Take any embedding of $G$ into the plane, which respects the orientation of edges around the nodes (edges may cross!). Then a subgraph $P$ is represented as a region  encircled with a dashed closed curve, drawn over the embedding of $G$ (namely $P$ is the collection of edges and nodes which are inside the bounded region delimited by the dashed curve), such that no node of $G$ is crossed by the dashed curve. 

\paragraph{Local and global moves.} On  $GRAPH$ we define moves (transformations) of a subgraph into another subgraph. Using the convention of drawing subgraphs, any such 
transformation is applied  only  inside the dashed curve ("inside" meaning the collection of bounded connected parts of the plane minus  the dashed curve).

Notice that the moves have to be defined such that they are independent on the embeddings in the plane. This is realized by numbering:  the edges which cross the curve (thus connecting the subgraph $P$ with the rest of the graph) will be numbered clockwise. Any edge which crosses the dashed curve which delimits the subgraph will bear one or several decorations 
(numbers) which will be used for gluing back the transformed subgraph with the remaining part of the graph (subjected to the move): after the transformation is performed, the edges of the transformed graph will connect to the graph outside the dashed curve by respecting the numbering of the edges which cross the dashed line.

These moves are local or global, according to the following criterion for "local".

For any   natural  number $N$ and any graph $G$ in $GRAPH$, let  
$\displaystyle \mathcal{P}(G,N)$ be the collection  of subgraphs $P$ of the graph $G$ which have the sum of the number of  edges and nodes less than or equal to $N$.

\begin{definition}
A local move has the following form: there is a number $N$ and a condition $C$ which is formulated in terms of subgraphs which have the sum of the number of  edges and nodes less than or equal to $N$,  such that for any graph $G$ in $GRAPH$ and for any $P \in \mathcal{P}(G,N)$, if $C$ is true for $P$ then transform $P$ into $P'$, where $P'$ is also a graph which have the sum of the number of  edges and nodes less than or equal to $N$.

Likewise, a subset of $GRAPH$ is defined by a local condition if it contains all graphs in $GRAPH$ which satisfy a condition which is  formulated in terms of subgraphs which have the sum of the number of  edges and nodes less than or equal to  a natural number $N$. 
\end{definition}

Moves or conditions  which are not local will be called "global". 

\paragraph{2.1. Graphic $\beta$ move.} This is the most powerful move in this graphic calculus, inspired by the $\beta$-reduction from untyped lambda calculus.

\centerline{\includegraphics[width=80mm]{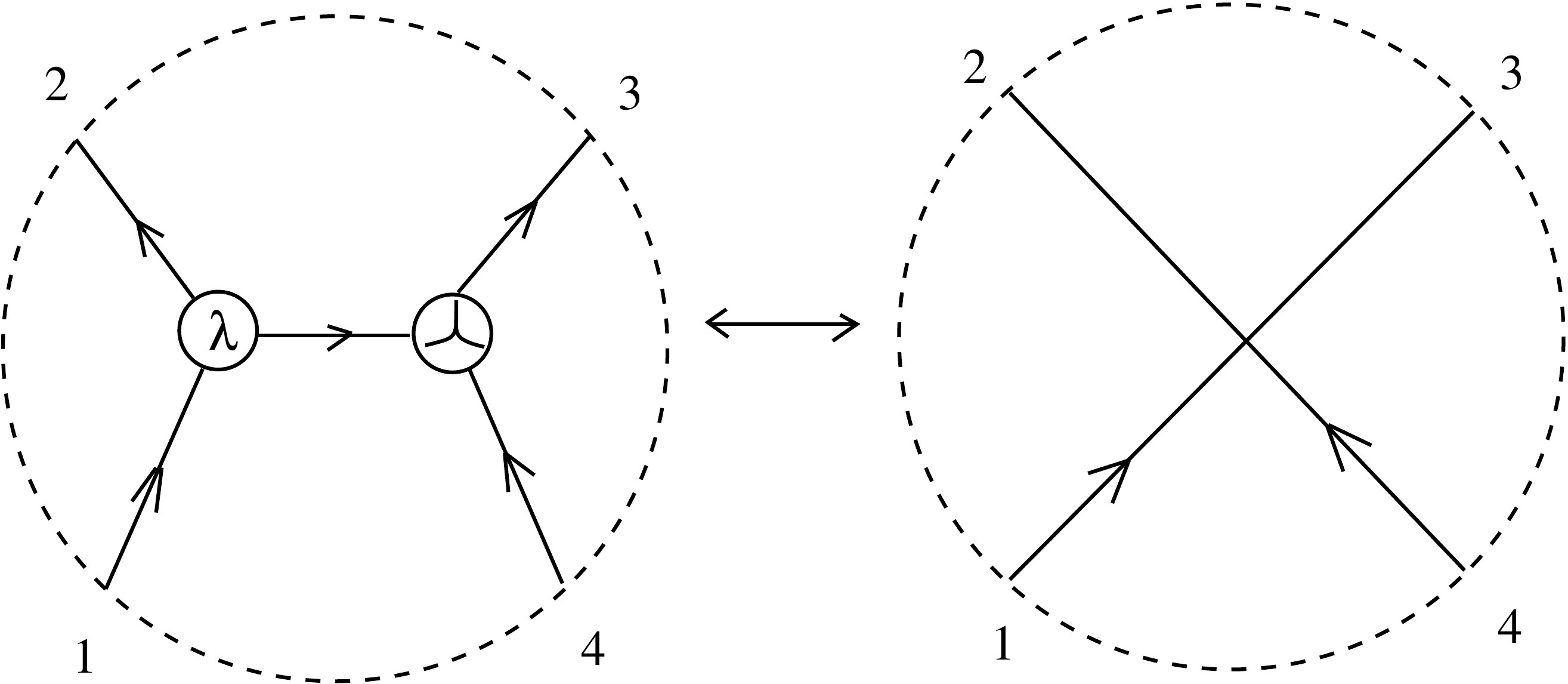}}

\begin{remark} There is a strong resemblance between this move and the "unzip operation" from the knotted trivalent graphs (KTGs) formalism (see \cite{ktg} section 3 and references therein). There are two differences: the unzip operation works on unoriented trivalent graphs and moreover by this move we can only unzip pairs of connected nodes (so the move goes only in one direction).
\end{remark}

\paragraph{2.2. The other local moves.} These are described in \cite{graphic} section 2: (CO-ASSOC), (CO-COMM), (R1), (R2), (ext2), local pruning. 

\paragraph{2.3. Global moves.} Described in the same section, these are: (ext1), (Global FAN-OUT), global pruning, elimination of loops.

\paragraph{2.4. Degenerate graphic $\beta$ moves obtained by elimination of loops.} There are  two situations where  the application of  the graphic $\beta$ move results in loops with no nodes.  In these cases we get graphs not in $GRAPH$, therefore we need to modify the $\beta$ move in order to take care of this.

\vspace{.5cm}

\centerline{\includegraphics[width=80mm]{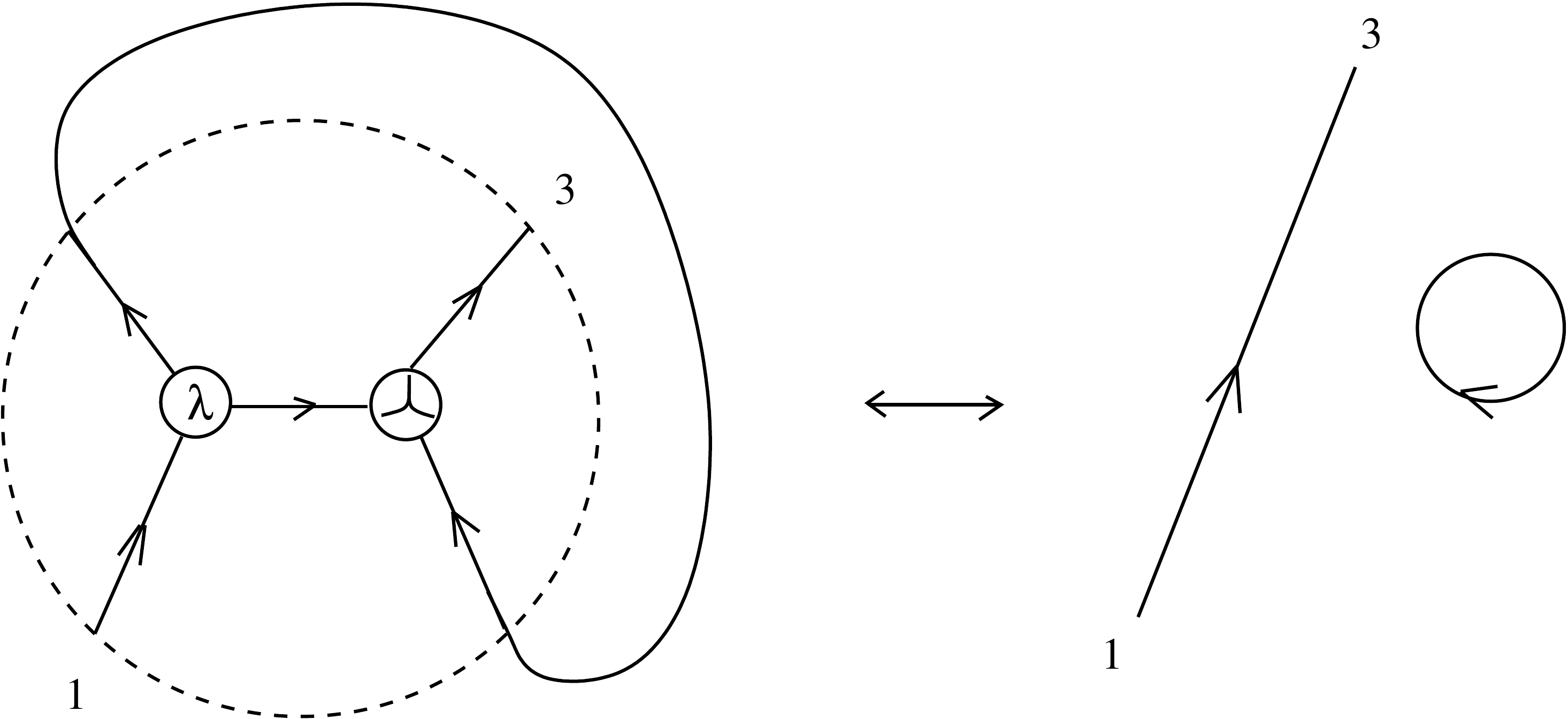}}

\vspace{.5cm}

\vspace{.5cm}

\centerline{\includegraphics[width=80mm]{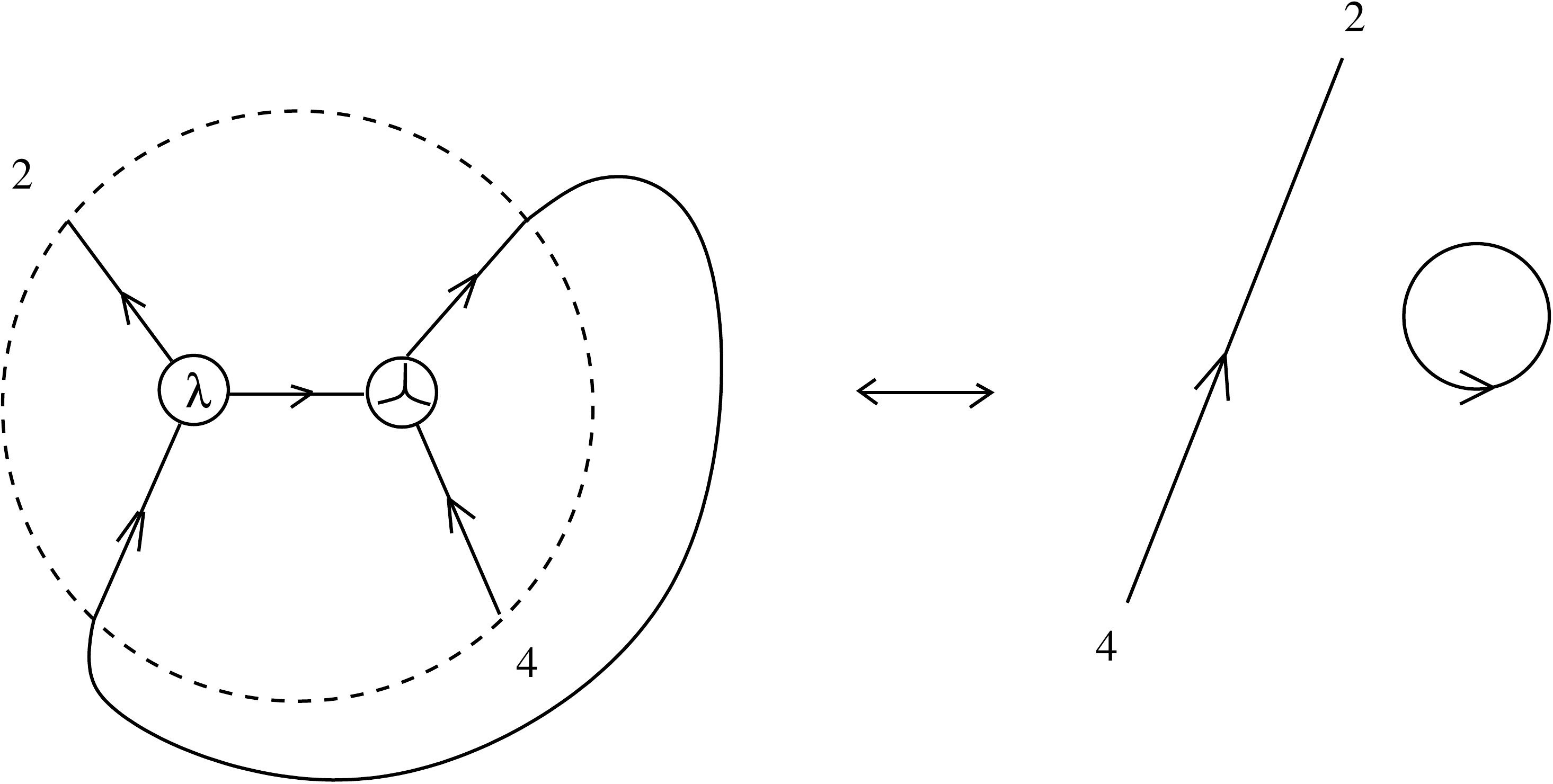}}

\vspace{.5cm}

 In these cases we get graphs not in $GRAPH$, therefore we need to modify the $\beta$ move in order to take care of this.  We do this by eliminating the loops.

The elimination of loops  takes the following form. 

\vspace{.5cm}

\centerline{\includegraphics[width=100mm]{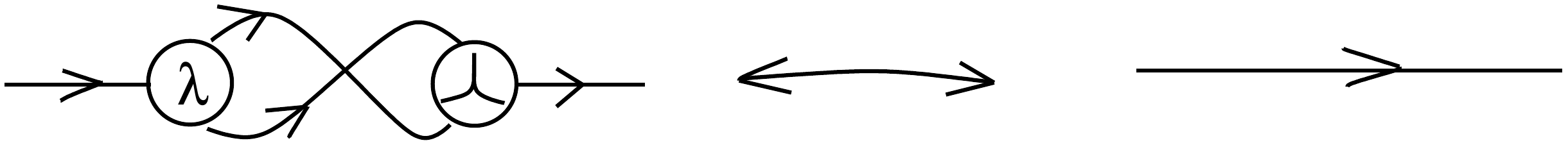}}

\vspace{.5cm}

\vspace{.5cm}

\centerline{\includegraphics[width=100mm]{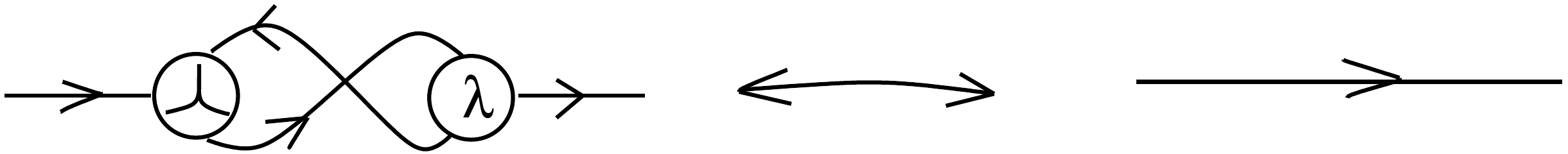}}

\vspace{.5cm}

Notice that there are other two cases of $\beta$ moves which don't pose any problem. 

\vspace{.5cm}

\centerline{\includegraphics[width=100mm]{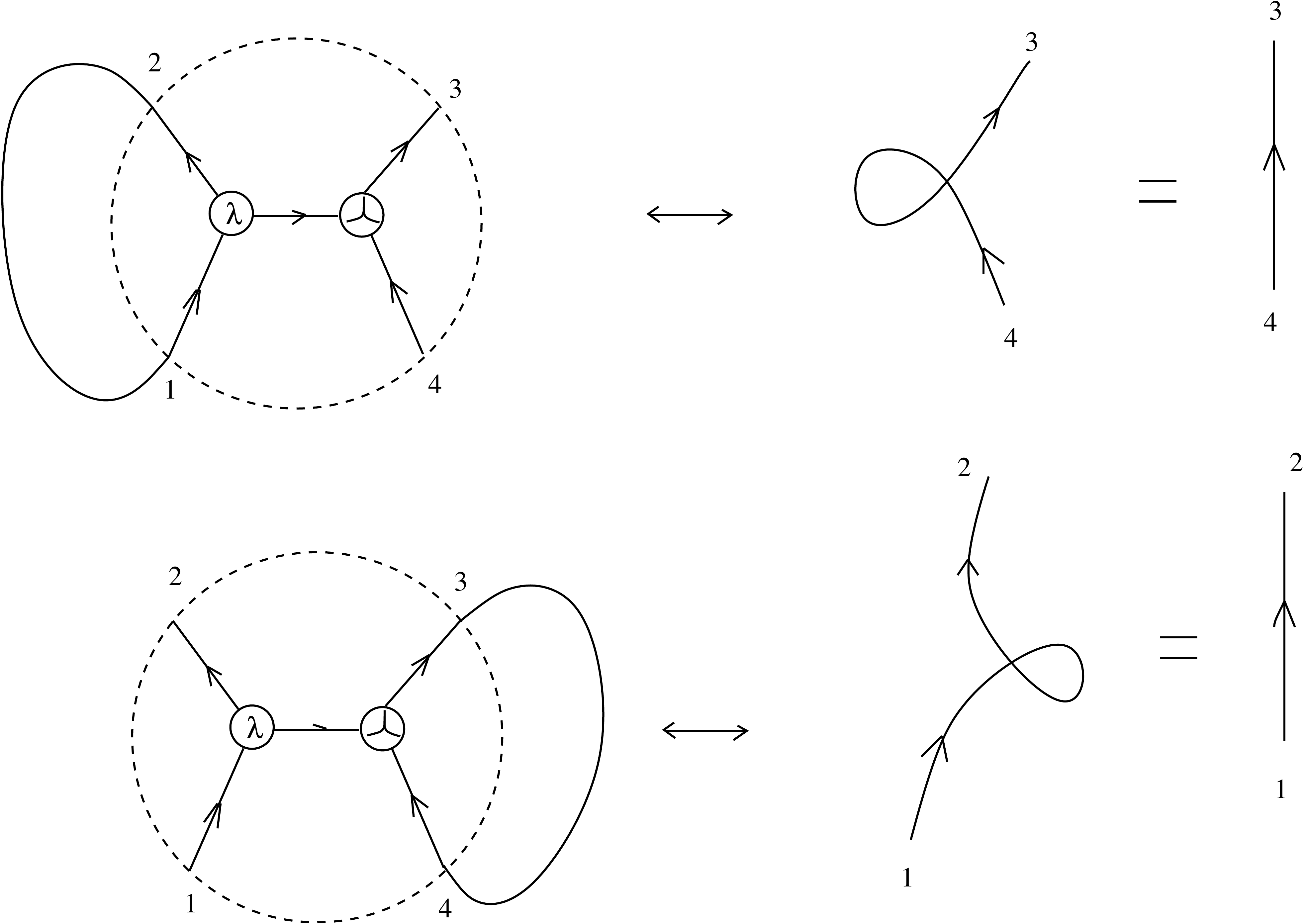}}

\vspace{.5cm}

\section{The lambda calculus sector} 

Untyped lambda calculus appears as a collection of moves applied to a subset of $GRAPH$ defined by a global condition. 

Let us number clockwise, by 1, 2, 3, the edges of the elementary graph $\lambda$  such that 1 is the number of the entrant edge. A $\lambda$-graph is a graph which belongs to the set $\lambda GRAPH$.   A $\lambda$-graph $G$ does not contain any  $\bar{\varepsilon}$ gate and moreover it satisfies the following global condition:

\begin{enumerate} 
\item[-] for any node $\lambda$ of the graph  any oriented path in $G$ starting at the edge 2 of this node can be completed to a path which either terminates in a graph $\top$, or else terminates at the edge 1 of the $\lambda$ node. 
\end{enumerate}

\begin{theorem}
(simplified version of theorem 3.1 \cite{graphic}) There is a transformation of terms in untyped lambda calculus into elements of the set $\lambda GRAPH$, which is bijective modulo  (CO-ASSOC), (CO-COMM), global pruning and (Global FAN-OUT) such that $\beta$ reduction is transformed into the  graphic $\beta$ move. Moreover, under this transformation , $\eta$ reduction becomes the (ext1) move. 
\end{theorem}

\section{The knot diagrams sector}
\label{secbraid}

In this section we use the graphic $\beta$ move in ways which are exterior to lambda calculus. 

Let us look again at the graphic beta move. In fact, this is a move which transforms a pair of edges (from the right of the picture) into the graph from the left of the picture. The fact that the edges cross in the figure is irrelevant. What matters is that, for the sake of performing the move, one edge is temporarily decorated with 1-3 and the other with 4-2.

Here are two more equivalent depictions of the same rule, coming from different choices of  1-3, 4-2 decorations. We may see the graphic beta move as: 

\vspace{.5cm}

\centerline{\includegraphics[width=80mm]{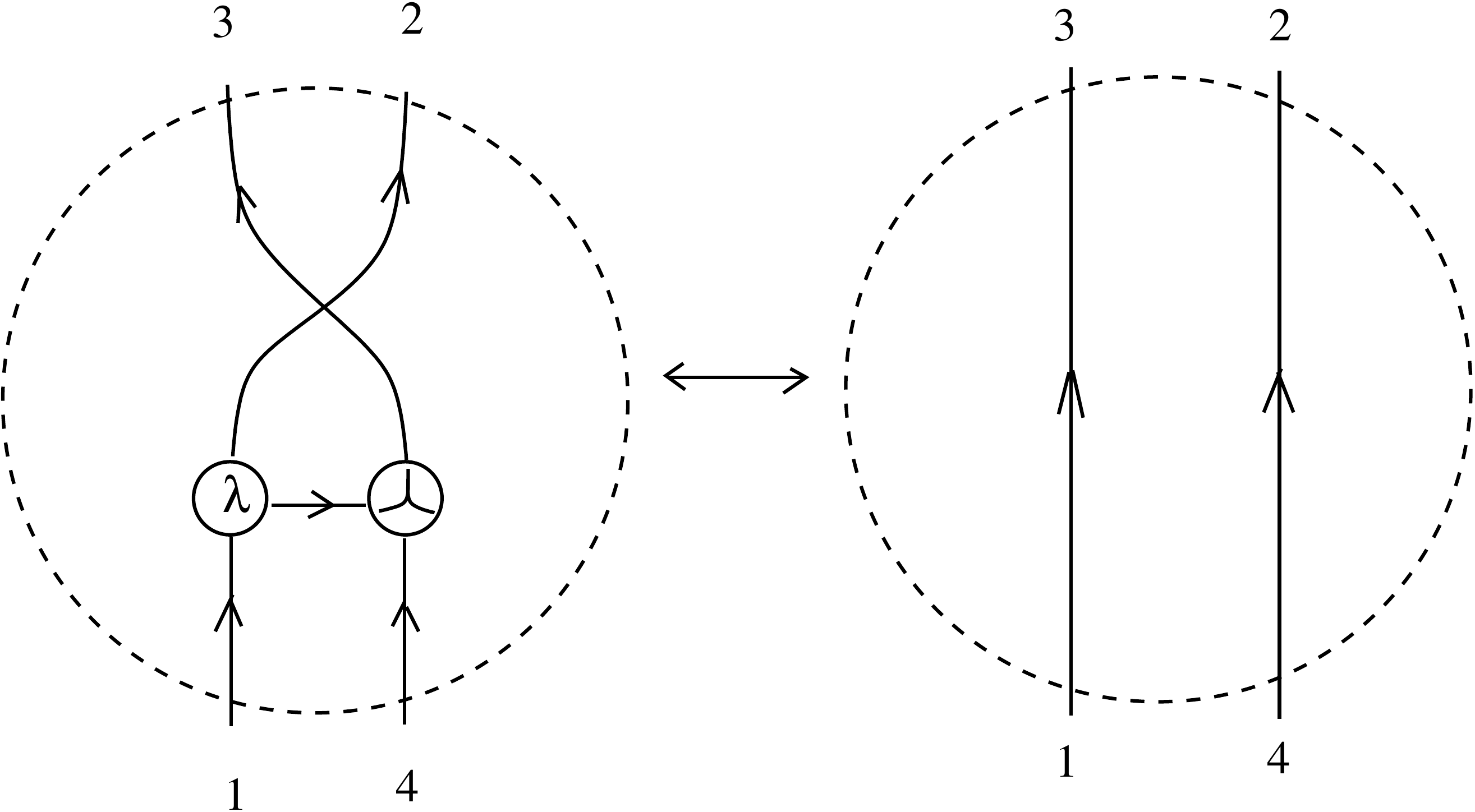}}

\vspace{.5cm}

By a different choice of the 1-3, 4-2 decoration, the graphic beta move can be also seen as: 

\vspace{.5cm}

\centerline{\includegraphics[width=80mm]{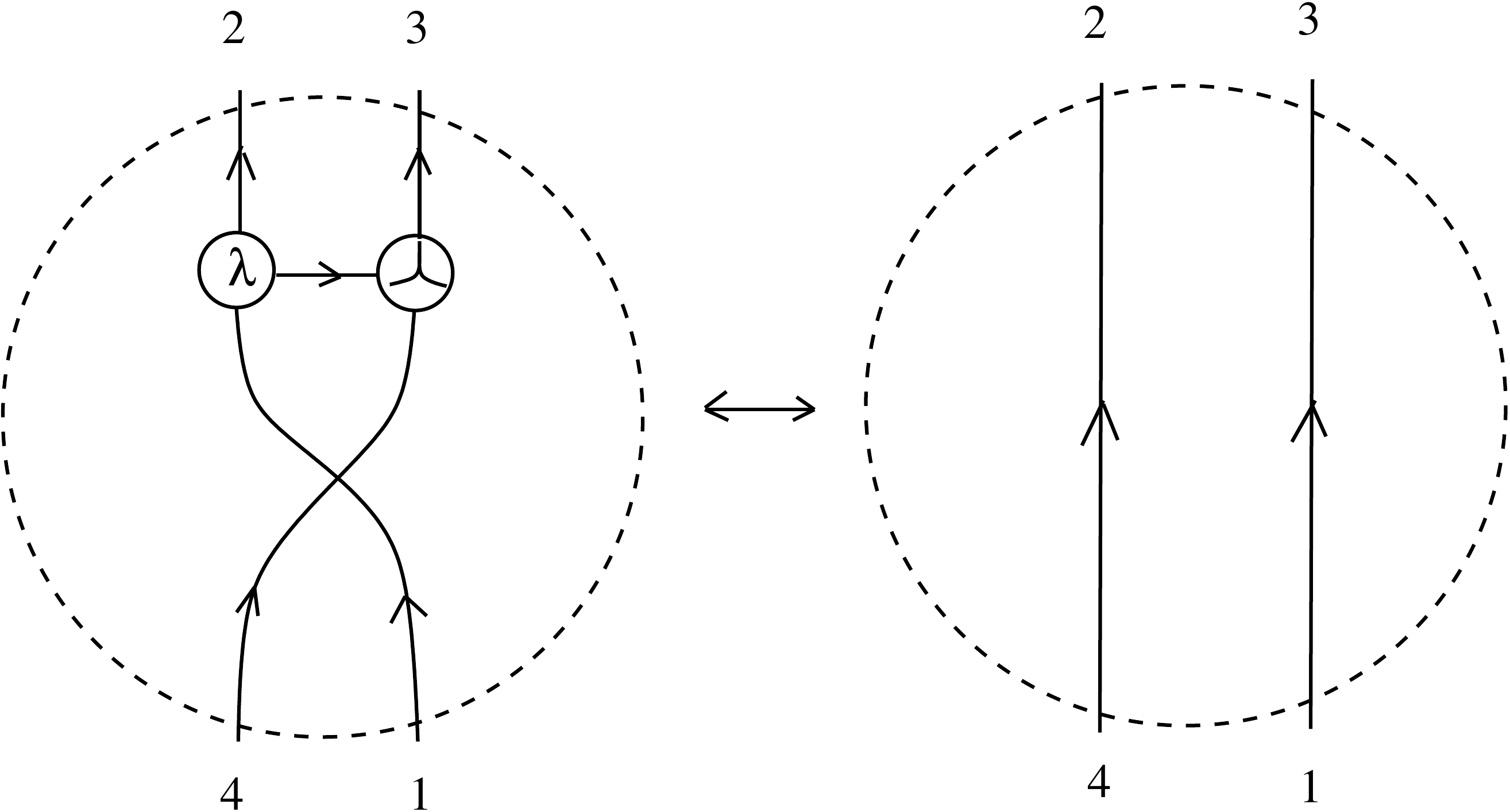}}

\vspace{.5cm}

Let us make then  some notations of the figures from the left hand sides of the previous diagrams.   Recall that a  tangle is a 1-dimensional oriented manifold embedded in the 3-dimensional space, up to regular isotopy. A link is an embedding of a finite collection of oriented circles in the 3-dimensional space  and a knot is an embedding of one oriented circle in 3D space (up to regular isotopy). Informally, links are knotted knots and tangles are like links, only that they may have loose ends. 

Diagrams of tangles are obtained via a regular projection on a 2D plane, together with additional over- and under- information at crossings.

\begin{definition}
We introduce the following two oriented crossings notations for the configuration from the left hand sides of the previous pictures: 

\vspace{.5cm}

\centerline{\includegraphics[width=80mm]{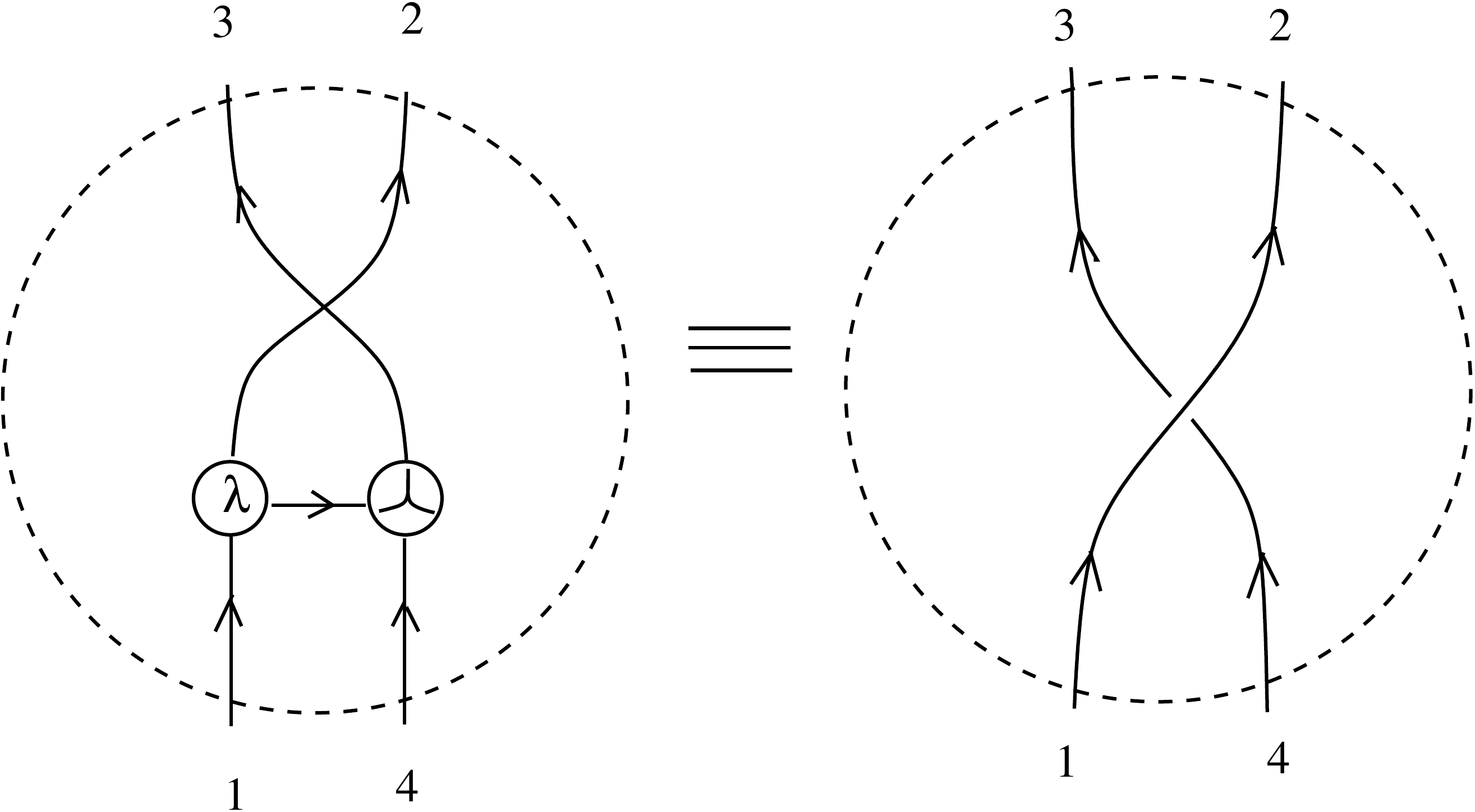}}

\vspace{.5cm}

\centerline{\includegraphics[width=80mm]{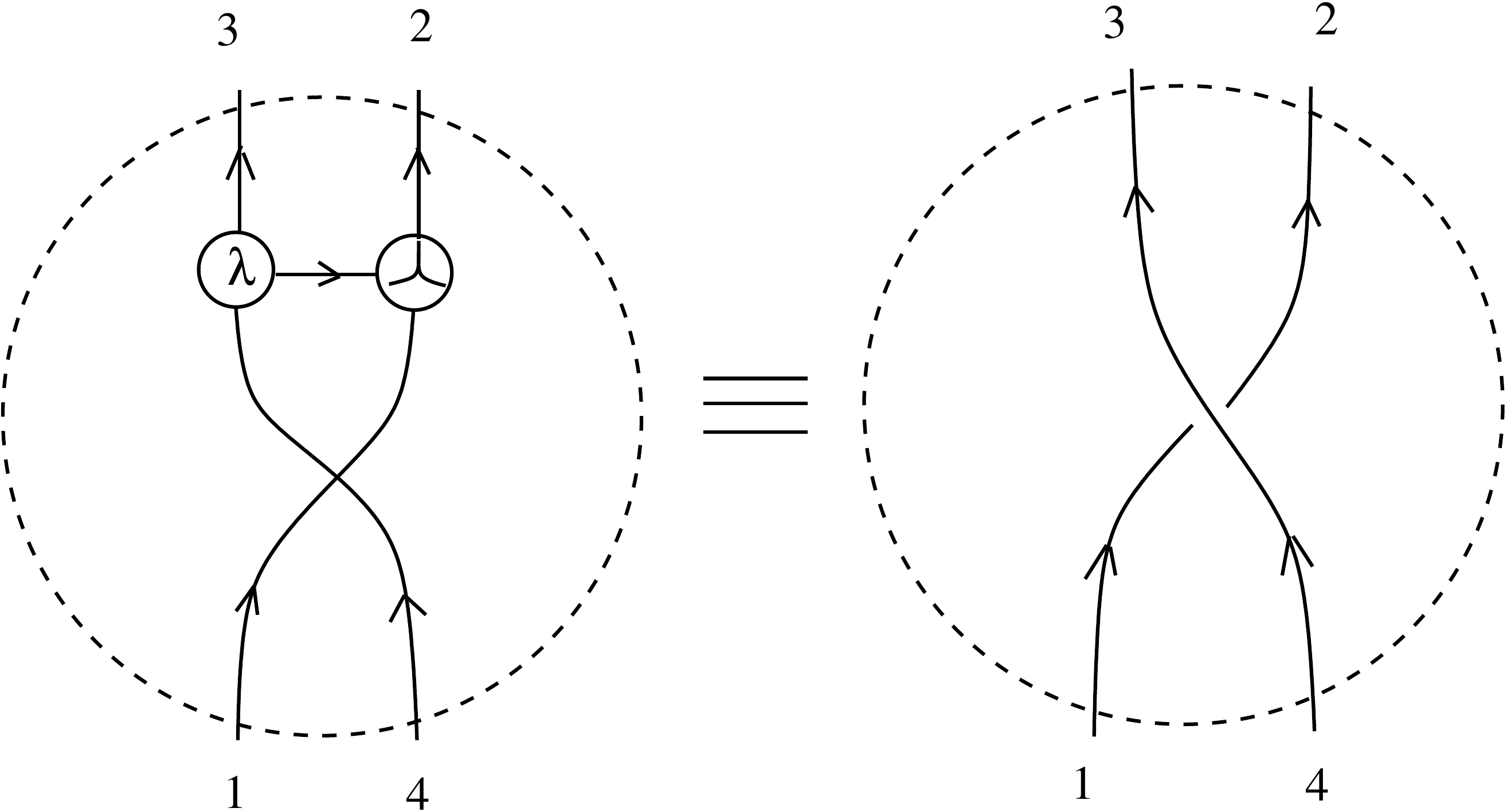}}

\vspace{.5cm}

The oriented crossings notations should be seen as "macros" in computer science sense, i.e. as a collection of "instructions" represented in an abbreviated format.

By using these notations, we define the subsets $TANGLEGRAPH$ and $LINKGRAPH$ of $GRAPH$ as the ones formed by all graphs in $GRAPH$ obtained from oriented tangle diagrams and  link diagrams respectively. 

Notice that $TANGLEGRAPH$ is defined by local conditions and $LINKGRAPH$ by global conditions. 

\label{crossingmacros}
\end{definition}

By using definition \ref{crossingmacros}, we may identify tangle diagrams (or link diagrams) with elements of $TANGLEGRAPH$ (or $LINKGRAPH$ respectively), called "tangle graphs". This vaguely specified "additional over- and under- information at crossings" takes a precise form by marking the string which passes "over" with a $\lambda$ (gate) and the string which passes "under" with a $\curlywedge$ (gate). Indeed, both crossing notations from definition \ref{crossingmacros} appear in 3D (i.e. after identifying tangle diagrams with elements of $TANGLEGRAPH$ and then interpreting tangle diagrams as projections of tangles) like this: 

\vspace{.5cm}

\centerline{\includegraphics[width=80mm]{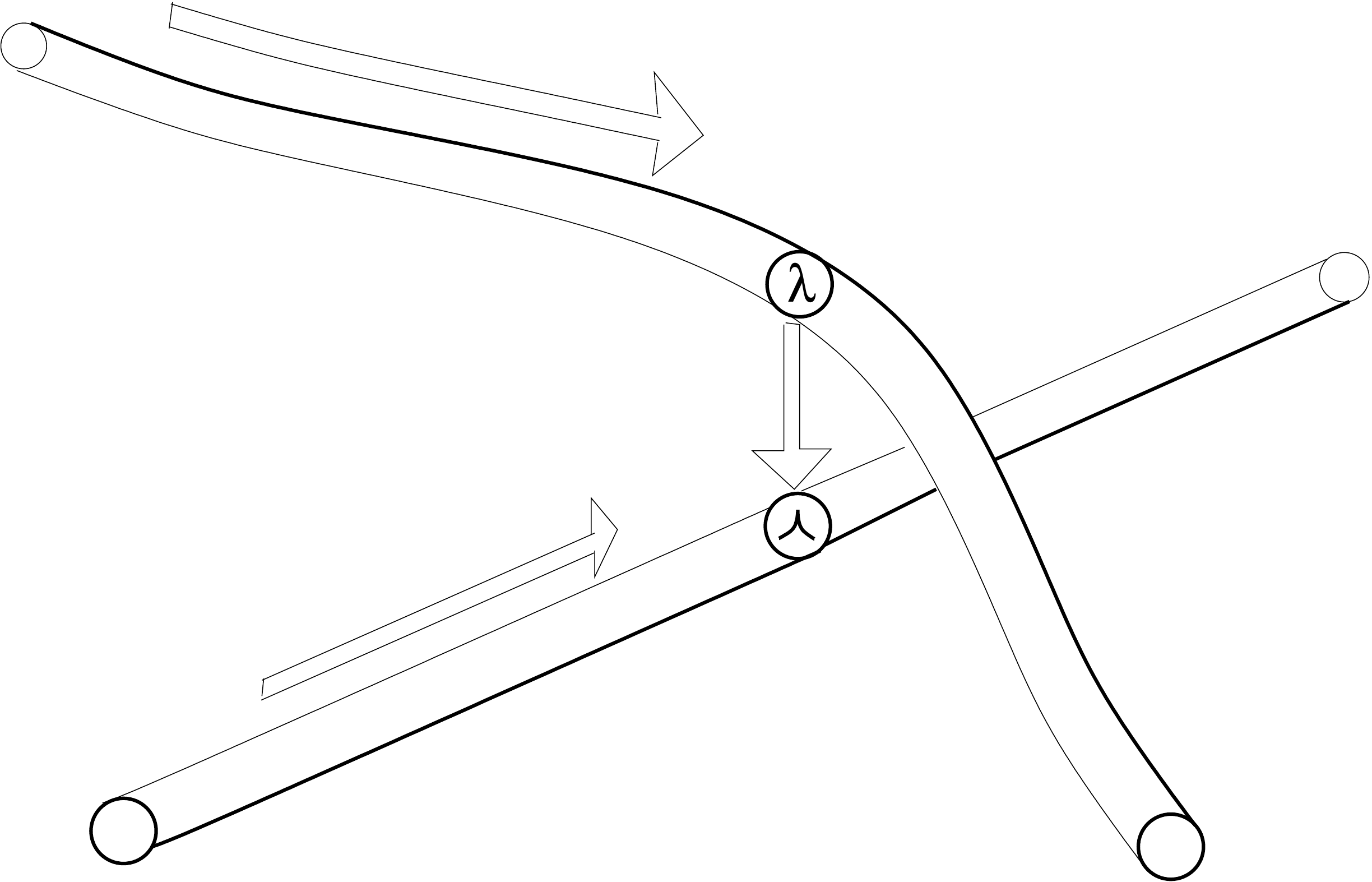}}

\vspace{.5cm}

The set $TANGLEGRAPH$ is defined by local conditions because any graph in the set  $TANGLEGRAPH$ is obtained from a tangle diagram. Tangle diagrams are locally defined objects, because there is no condition on the "loose ends". On the contrary, a link diagram is a tangle diagram which satisfies the global condition of non-existence of loose ends. That is why $LINKGRAPH$ is a set defined by global conditions.  

\begin{theorem}
(a) The set $TANGLEGRAPH$ is closed under graphic $\beta$ move and elimination of loops.

(b) Two link diagrams represent projections of the same link (up to regular isotopy) if and only if we can pass from the associated graph graph of the first diagram to the associated link graph of the second diagram by using a finite number of local moves from the following list: 
\begin{enumerate}
\item[-] the two elimination of loops moves, call them (R1a) 
\vspace{.5cm}

\centerline{\includegraphics[width=100mm]{betar_r1a.pdf}}

\vspace{.5cm}

and (R1b) respectively,

\vspace{.5cm}

\centerline{\includegraphics[width=100mm]{betar_r1b.pdf}}

\vspace{.5cm}

\item[-] the move called (R2a), which is a composite of two graphic $\beta$ moves 

\vspace{.5cm}

\centerline{\includegraphics[width=100mm]{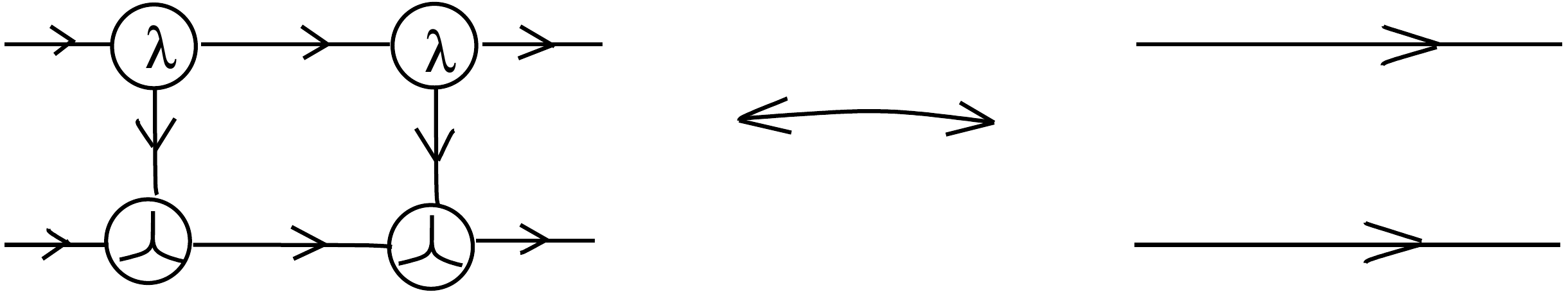}}

\vspace{.5cm}

\item[-] and the move called (R3a), which is a composite of 6 graphic $\beta$ moves 

 \vspace{.5cm}

\centerline{\includegraphics[width=100mm]{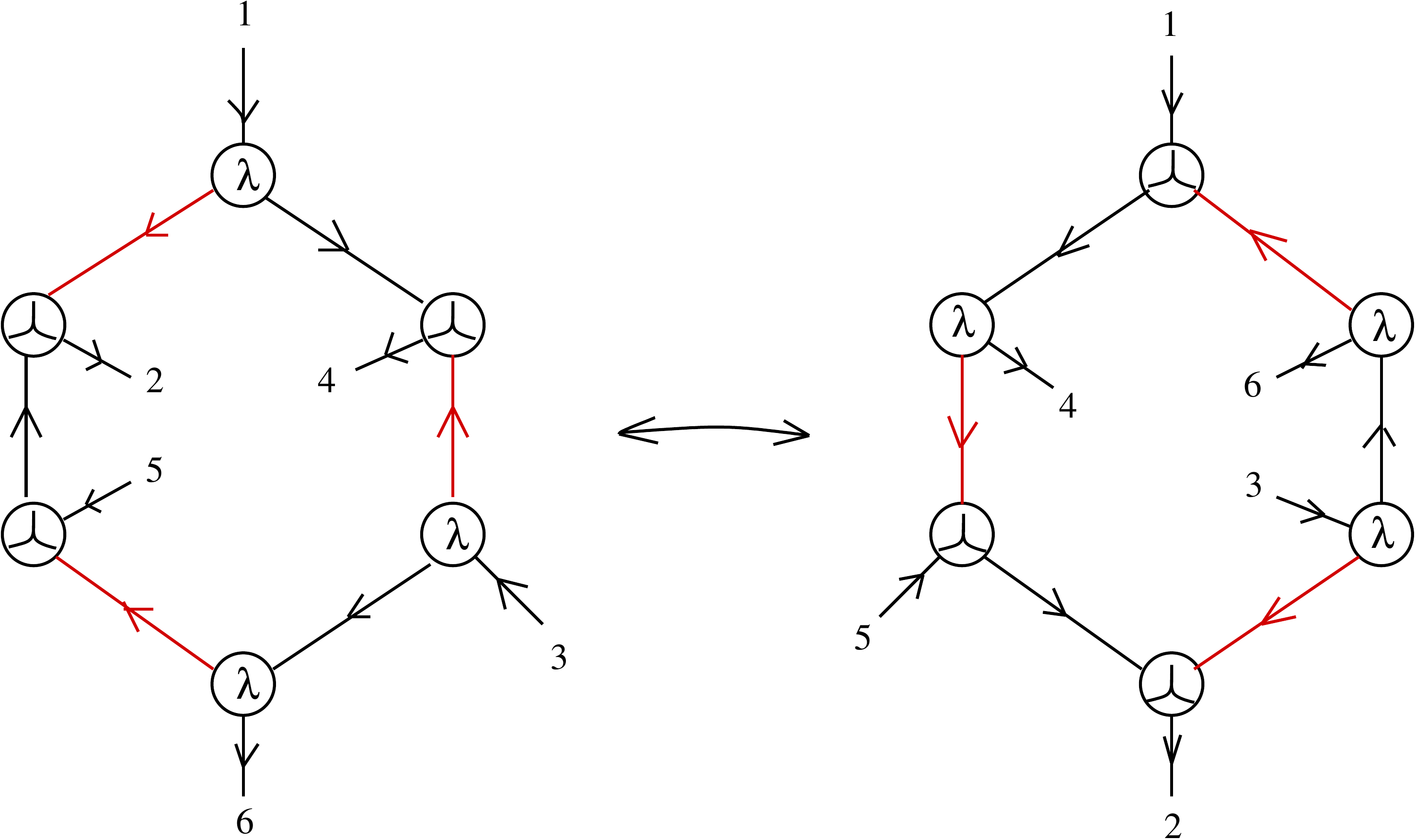}}

\vspace{.5cm}
(the red arrows show which pairs of gates $\lambda$ - $\curlywedge$ give the over-under information for the associated crossings, see further). 
\end{enumerate}

\end{theorem}

\paragraph{Proof.} (a) With the notations from definition \ref{crossingmacros}, the graphic 
$\beta$ move appears under the following two forms: 

\vspace{.5cm}

\centerline{\includegraphics[width=100mm]{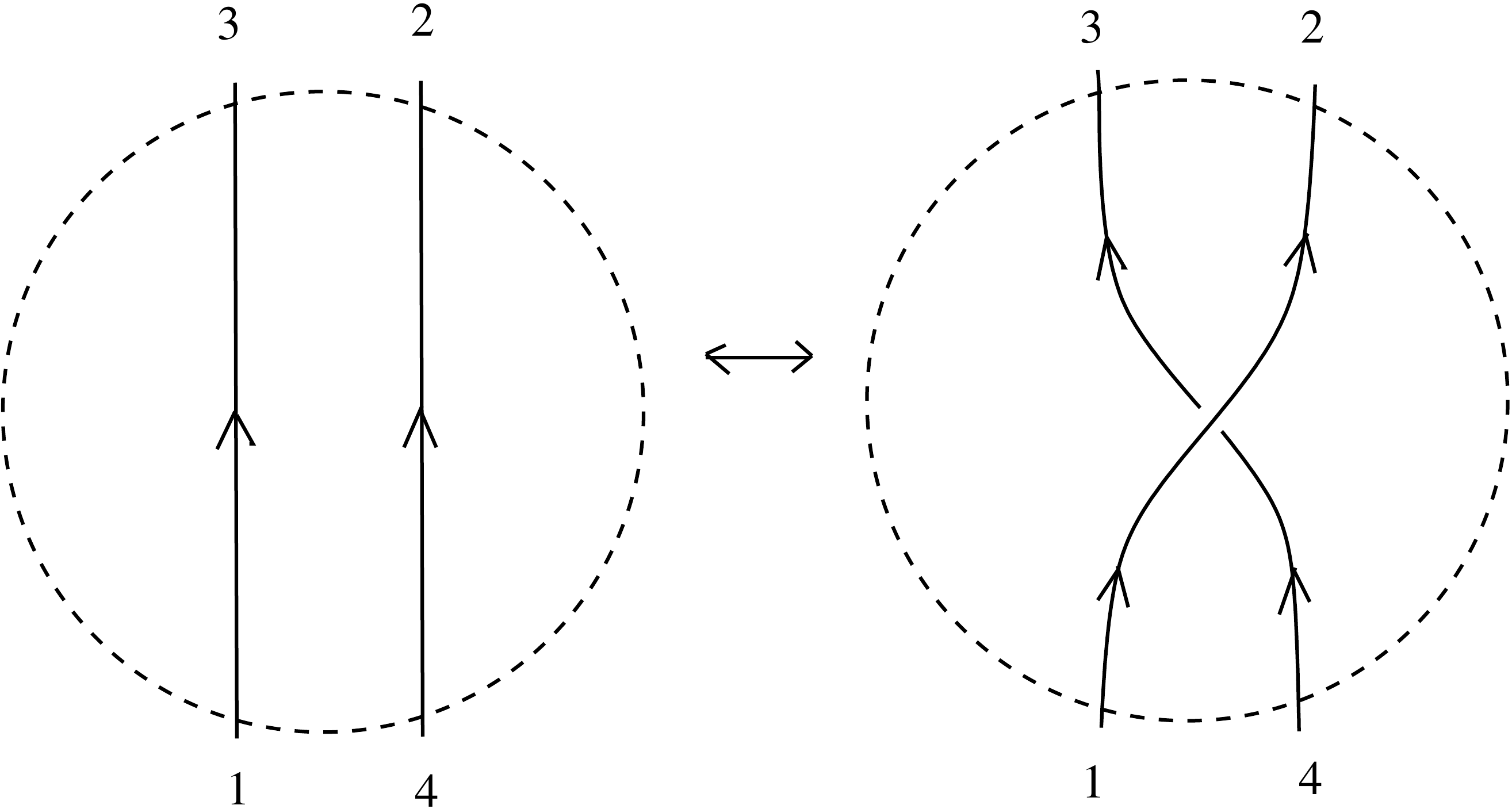}}

\vspace{.5cm}

\centerline{\includegraphics[width=100mm]{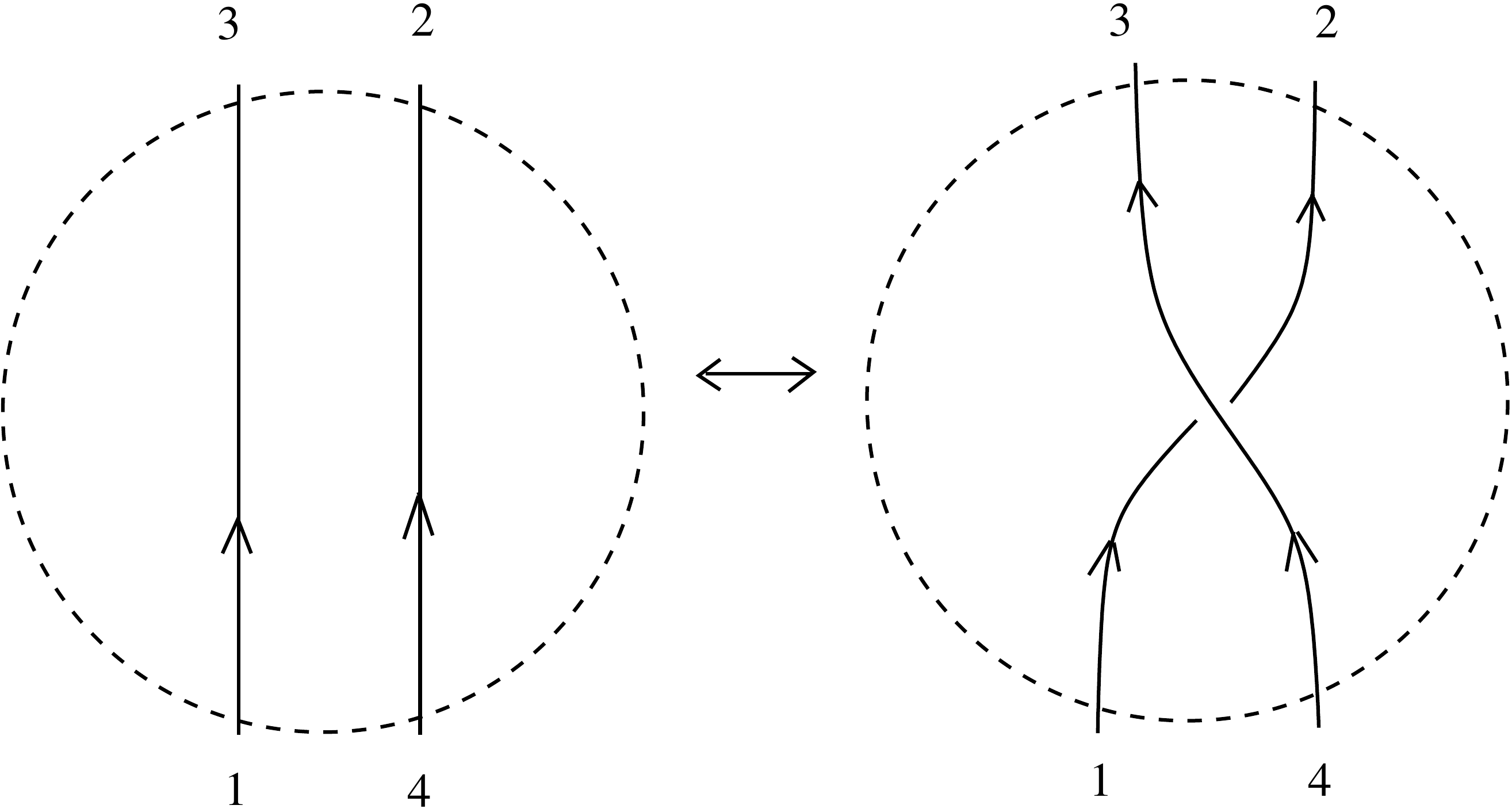}}

\vspace{.5cm}
therefore the graphic $\beta$ move appears as a braiding or an unbraiding operation, in the case of tangle graphs. In the case of link graphs, they may be transformed into  tangle graphs. Notice also that the only possibility of applying the graphic $\beta$ move is in this context, because tangle graphs contain arrows from the right exit peg of a gate $\lambda$ to the left input peg of a gate $\curlywedge$ only when such a pair encodes the over-under relation of a crossing. Therefore the only possible graphic $\beta$ moves are those which are related to braiding or unbraiding a crossing. 

For the closure with respect to the elimination of moves, see point (b), where it is shown that these moves correspond to Reidemeister I move(s) for oriented crossings. We use also the observation previously made in order to exclude any other instance of application of elimination of loops than the ones related to the Reidemeister I moves. 

(b) Two link diagrams represent the same link if and only if one can be obtained from the other by using a finite number of Reidemeister moves.   We shall follow Polyak \cite{polyak}, with the only modification of notation consisting in using a capital "R" instead of "$\Omega$" for the names of Reidemeister moves.  For oriented link diagrams, there are 24 Reidemeister moves, but a minimal generating set is formed by the moves (R1a), (R1b), (R2a), (R3a), according to theorem 1.1 \cite{polyak}. Our purpose is to show that the Reidemeister moves (R1a), (R1b), (R2a), (R3a) correspond to the local moves in $LINKGRAPH$  with the same names which are listed in the text of the theorem. 

Further we give the concerned Reidemeister moves and we leave to the interested reader to make the translation in terms of link graphs, by using definition \ref{crossingmacros}. (When needed, we shall number the input and output leaves, following the same convention as the one used for graphs in $GRAPH$.)

The Reidemeister move (R1a) is the following: 

\vspace{.5cm}

\centerline{\includegraphics[width=100mm]{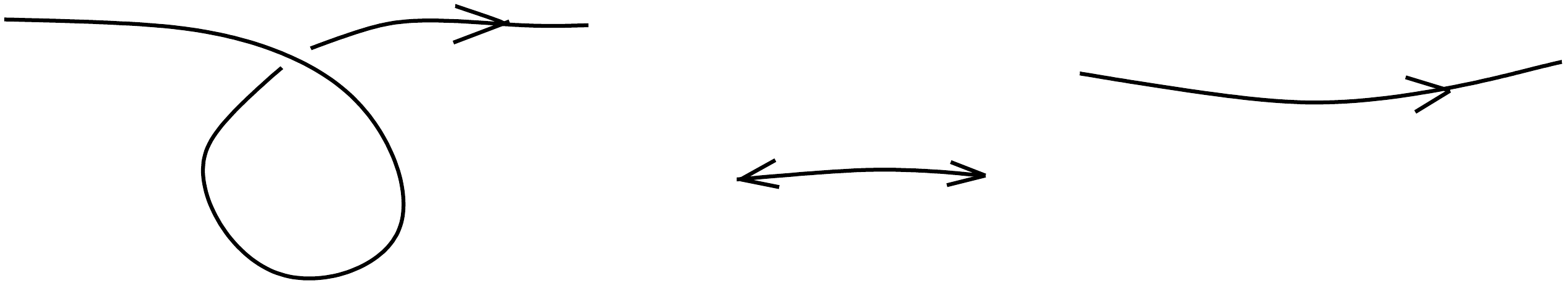}}

\vspace{.5cm}
which translates into the elimination of loops move (R1a) in $LINKGRAPH$. In the same way, the Reidemeister move (R1b) is 

\vspace{.5cm}

\centerline{\includegraphics[width=100mm]{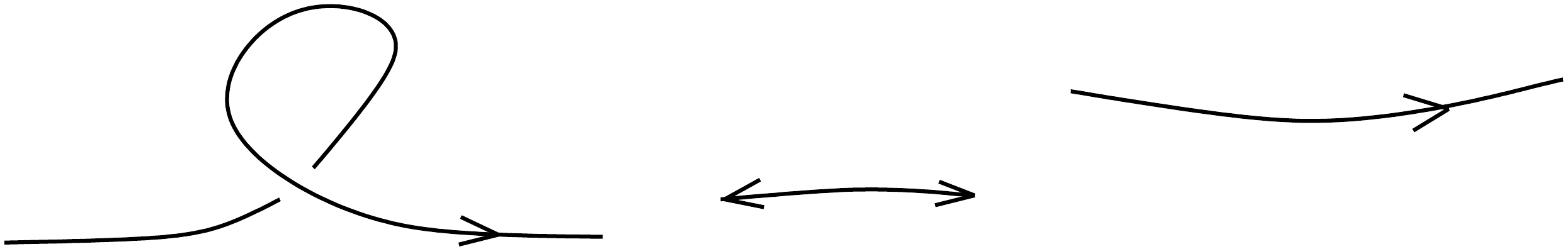}}

\vspace{.5cm} 
and it becomes the elimination of loops move (R1b). The Reidemeister move (R2a) 

\vspace{.5cm}

\centerline{\includegraphics[width=100mm]{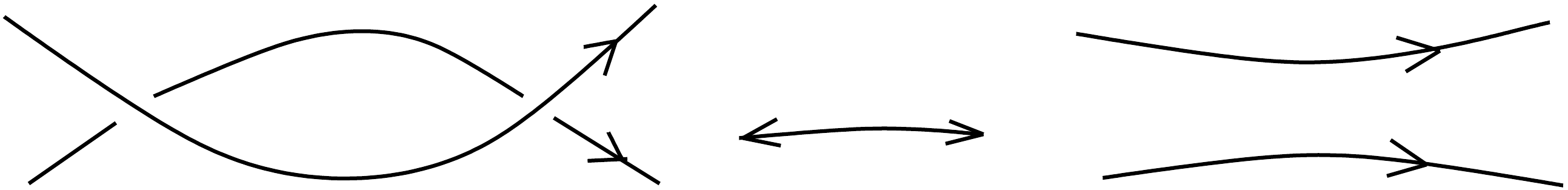}}

\vspace{.5cm}
becomes the move (R2a), which can be seen as a composition of two graphic $\beta$ moves, in obvious way. 

Finally, the Reidemeister move (R3a) 

\vspace{.5cm}

\centerline{\includegraphics[width=100mm]{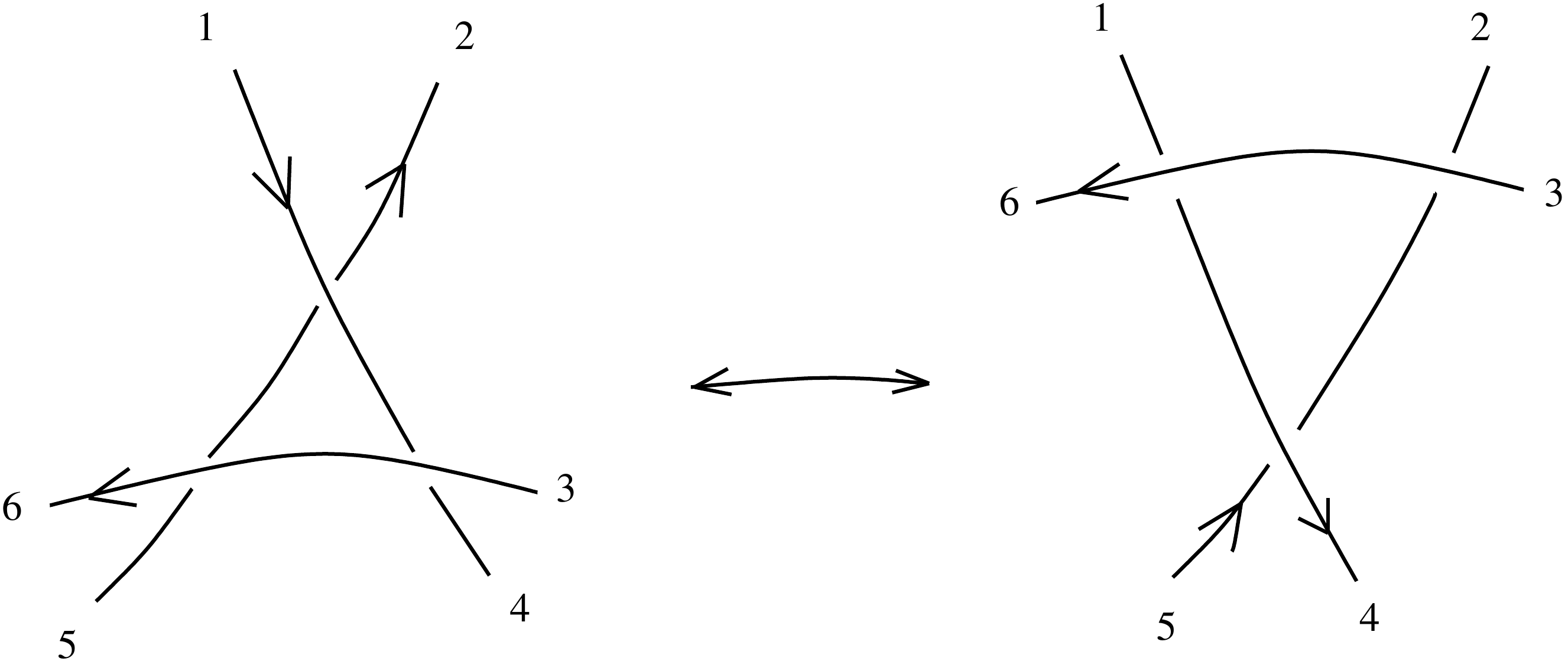}}

\vspace{.5cm}
becomes the move (R3a) from the text of the theorem. This move can be seen as a composition of 6 $\beta$ moves, in several ways. \quad $\square$

\section{Conclusion} 

In this paper we put on common ground untyped lambda calculus and link diagrams. For a very interesting discussion about lambda calculus and link diagrams see Kauffman \cite{kauf}, especially pages 46-48. On a related note, Meredith and Snyder \cite{mersny} associate to any knot diagram a process in the $\pi$-calculus.

To any oriented link diagram we associated a graph in $GRAPH$, such that the Reidemeister moves translate into particular compositions (of an even number) of graphic $\beta$ moves and any number of elimination of loops moves. The family of these graphs, called link graphs, form a sector of the graphic lambda calculus, different than the sector of the same calculus which covers untyped lambda calculus.

\end{document}